\documentclass[11pt,oneside]{article}
\usepackage[utf8]{inputenc}
\usepackage[english]{babel}
\usepackage{amsmath}
\usepackage{amsfonts}
\usepackage{amssymb}
\usepackage{amsthm}
\usepackage{tikz-cd} 
\usepackage{hyperref}
\usepackage{pgfplots}
\usetikzlibrary{calc}

\usepackage[left=3cm,right=3cm,top=3cm,bottom=3cm]{geometry}
\usepackage[title]{appendix}
\usepackage[]{algorithm2e}

\usepackage[normalem]{ulem}
\useunder{\uline}{\ul}{}
\theoremstyle{plain}
\newtheorem{teo}{}[section]
\newtheorem{prop}[teo]{Proposition}
\newtheorem{cor}[teo]{Corollary}

\newtheorem{lem}[teo]{Lemma}
\newtheorem{thm}[teo]{Theorem}
\theoremstyle{definition}
\newtheorem{ex}[teo]{Example}
\newtheorem{rem}[teo]{Remark}
\newtheorem{df}[teo]{Definition}

\newcommand\blfootnote[1]{%
  \begingroup
  \renewcommand\thefootnote{}\footnote{#1}%
  \addtocounter{footnote}{-1}%
  \endgroup
}

\definecolor{cof}{RGB}{219,144,71}
\definecolor{pur}{RGB}{186,146,162}
\definecolor{greeo}{RGB}{91,173,69}
\definecolor{greet}{RGB}{52,111,72}

\usepackage{graphicx}

\title{Negative-Dimensional Simplicial Complexes via the Grothendieck Group of the Join Monoid}
\author{Pedro J. Chocano and L. Felipe Prieto-Martínez}
\date{}

\begin{document}
\maketitle

\begin{abstract}
We introduce a framework that extends classical combinatorial and topological notions of simplicial complexes to a setting that includes simplicial complexes of negative dimension. This is achieved by considering the Grothendieck group of the commutative monoid of finite simplicial complexes under the join operation, which provides formal inverses and allows negative-dimensional objects to be defined. Within this group we adapt standard invariants such as $f$‑ and $h$‑vectors, their corresponding polynomials, and the Euler characteristic.
The same construction naturally endows the Grothendieck group with an Alexandroff topology. We show that this space is universal for finite topological spaces, and consequently the associated ordered simplicial complex is universal within its class. 
\end{abstract}
\blfootnote{This research is partially supported by Grant  PID2024-156663NB-I00 from Ministerio de Ciencia, Innovación y Universidades (Spain) and  2025/SOLCON-159637 from Rey Juan Carlos University}

\section{Introduction}
The aim of this work is to develop a systematic framework in which classical notions from the theory of simplicial complexes—both combinatorial and topological—extend naturally to the setting of simplicial complexes of negative dimension. Our main observation is that the monoid of simplicial complexes under the join operation admits a canonical Grothendieck group. Within this group, formal inverses of joins acquire a geometric meaning: they behave as simplicial complexes whose dimensions are negative. This perspective allows us to reinterpret familiar invariants in a broader algebraic–combinatorial setting and to compute them for objects that do not exist in the classical theory.

Our starting point is the topological join~$*$. This construction assigns to two simplicial complexes $X$ and $Y$ a new simplicial complex $X * Y$ whose dimension is strictly larger than those of $X$ and $Y$. Standard examples include the suspension and the cone of a simplicial complex $X$, obtained respectively as $X * [2]$ and $X * [1]$, where $[2]$ denotes the discrete two‑point simplicial complex and $[1]$ the one‑point simplicial complex. Although fundamental in topology, the join operation lacks natural inverses: unless a simplicial complex is already known to decompose as $Z = X * Y$, there is no canonical procedure to “undo’’ the join. When such a decomposition is available, however, one may formally “desjoin’’ $Z$ with $Y$ to recover $X$.

The key insight is that the class of finite simplicial complexes equipped with the join satisfies the axioms required to form a commutative monoid. This makes it possible to apply the classical Grothendieck construction, which formally adjoins inverses with respect to the monoid operation. The resulting Abelian group provides a natural ambient space in which joins become invertible, and where one may interpret the formal inverse of a cone or a suspension as the join of a simplicial complex of negative dimension. Working in this setting, we show that several standard combinatorial and topological invariants—such as $f$‑vectors and $f$‑polynomials, $h$‑vectors and $h$‑polynomials, or the Euler characteristic—extend coherently to these new objects.

\begin{figure}[h]
    \centering
    \includegraphics[width=0.9\linewidth]{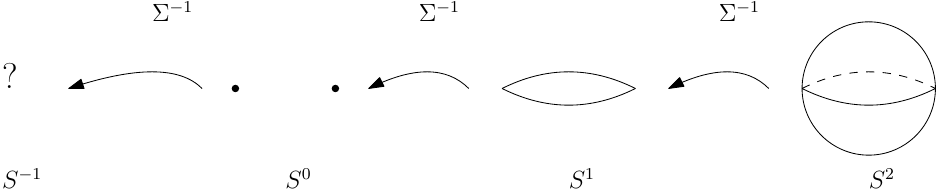}
    \caption{Schematic drawing of a “desuspension” from ~\cite{wolcott_mcternan_2012}.}
    \label{fig:idea_negativos}
\end{figure}

One of the motivations for this program came from the question “Desuspension: what does it lead to?” posed in~\cite{wolcott_mcternan_2012}, where the authors include a diagram closely resembling Figure~\ref{fig:idea_negativos}. In light of this question, our results may be understood as providing at least a partial answer: we are now able to compute and interpret invariants associated to negative simplicial complexes, such as the number of faces of the negative sphere $S^{-1}$ in each dimension. Nevertheless, a complete geometric understanding of these constructions remains open and is the subject of current research.

As an illustration, consider the negative dimensional sphere obtained by iterating the formal inverse of the suspension. We denote this space by $[\emptyset, S^n]$, where $n$ is a positive integer. Its dimension equals $-n-1$, and its total reduced Betti number is $\tfrac{1}{2}$. The extended $f$‑polynomial is
$$f([\emptyset,S^n])=\sum_{i=-\infty}^{-n-1}(-1)^{i+n+1}\binom{-i-1}{n}2^it^i,$$
and the corresponding extended $f$‑vector is
$$(\cdots,\underbrace{-\frac{(n+3)(n+2)(n+1)}{6\cdot2^{n+4}}}_{-n-4},\underbrace{\frac{(n+2)(n+1)}{2\cdot 2^{n+3}}}_{-n-3},\underbrace{-\frac{n+1}{2^{n+2}}}_{-n-2},\underbrace{\frac{1}{2^{n+1}}}_{-n-1},0,0,0,\cdots).$$
It follows that the Euler characteristic satisfies
$$\chi([\emptyset,S^n])=\begin{cases}
    0 & \textnormal{if } n \textnormal{ is odd,}\\
    2 & \textnormal{if } n \textnormal{ is even.}
\end{cases}$$

This example also suggests that techniques from the theory of Riordan matrices may be applied to obtain combinatorial information about the extended $f$‑vectors of families of spaces obtained through successive “desjoins’’. In particular, \cite{chocano} presents several specific families of spaces for which these techniques appear to be adaptable. Such methods could potentially shed light on recurrence relations, generating functions, and structural patterns arising in these extended invariants.

Although our original goal was to study simplicial complexes of negative dimension, the construction developed in this article has a further consequence: it induces a canonical Alexandroff topology on the Grothendieck group. Remarkably, the resulting space is universal for the class of finite topological spaces. Applying the associated functor shows that the simplicial complex naturally associated to the Grothendieck group is universal for ordered simplicial complexes. From this viewpoint, the Grothendieck group endowed with its Alexandroff topology may be regarded as a moduli space parametrizing both classical and negative simplicial complexes.

Beyond the results presented here, this work should be regarded as a first step in a broader program. The authors are currently developing, among other directions, a homology theory for negative simplicial complexes, as well as a more detailed study of several constructions introduced in this paper. These investigations aim to clarify the geometric meaning of negative-dimensional phenomena and to explore potential applications in combinatorial and topological settings.

The paper is organized as follows. Section~\ref{sec_preliminaries} gathers notation and standard definitions concerning monoids, simplicial complexes, Grothendieck groups, Alexandroff spaces, and other preliminary material, in order to keep the manuscript as self‑contained as possible. Section~\ref{sec:groth_group_join} studies the Grothendieck group associated with the monoid of finite simplicial complexes under the join operation. Section~\ref{sec_grothendieck_topology} introduces a natural Alexandroff topology on this group and investigates its principal properties; among them, we prove that the space constructed in the preceding section is universal for finite topological spaces. Finally, in Section~\ref{sec_homomorphism} we extend classical invariants of simplicial complexes to elements of the Grothendieck group, including, for example, the characteristic polynomial of a simplicial complex. This allows us to define simplicial complexes of negative dimension and to analyze combinatorial invariants such as their $f$‑vectors.

\section{Preliminaries}\label{sec_preliminaries}
\subsection{Monoids}\label{subsec_monoid_grothendieck}

We refer the reader to \cite{bourbaki, hungerford, jacobson, lang} for the standard concepts of monoid theory, but let us recall some of them, neccesary for the rest of this paper.

A non-empty set $M$, together with a binary operation $\cdot :M\times M\rightarrow M$, forms a \emph{monoid} if it possesses an identity element $e$ and the operation is associative. Such structure  is denoted by $(M,\cdot)$. Let $M,M'$ be monoids, a \emph{monoid homomorphism} is a mapping $f:M\rightarrow M'$ such that $f(a\cdot b)=f(a)\cdot f(b)$ for every $a,b\in M$. From now on, we will omit $\cdot$ in $a\cdot b$ and simply write $ab$. An element $a\in M$ is a \emph{unit} if it is invertible. A non-unit is called \emph{irreducible} if it can not be represented as a product of two non-units.
We say that:

\begin{itemize}
    \item The monoid is \emph{cancellative} if the \emph{cancellation law} holds, that is, whenever $a,b,c\in M$ and $a c=b c$, it follows that $a=b$. 
    \item The monoid is \emph{commutative} if the \emph{commutative law} holds, that is, for every $a,b\in M$, $ab=ba$.
    \item The monoid is a \emph{unique factorization monoid} if every non-unit can be expressed as a finite product of irreducible elements, and this factorization is unique up to the order of the factors and multiplication by units. 
    
    Every unique factorization monoid is cancellative.

    \item The monoid is \emph{positive} if the only unit is the identity element.

\end{itemize}

For unique factorization monoids, we define \emph{divisor} of some element in the usual way. Also, the concepts of \emph{irreducible element} and \emph{prime element} are equivalent, and we can also define the \emph{greatest common divisor} of two or more elements as expected. More about this will be said in the following section.

A monoid that is cancellative and positive is called a monoid \emph{with a natural order}, because of the following result:

\begin{lem} If $(M,\cdot )$ is a monoid that is cancellative and positive, the binary relation $\leq$ given by
\begin{equation} \label{eq.leq} a\leq b \qquad \Longleftrightarrow \qquad \exists d\in M,\; a d=b,\end{equation}

\noindent is an order relation that is compatible with the monoid operation, i.e., if $a\leq b$ and $c\in M$, then $ac\leq bc$. Moreover, the identity element is always a minimum.

\end{lem}

\begin{proof}
    (i) Reflexive property. Clearly $x\leq x$ since $x=x e$, where $e$ denotes the identity element. (ii) Antisymmetric property. Suppose $x\leq y$ and $y\leq x$, then there exists $d$ and $d'$ such that $y=x d$ and $x=y d'$. Hence, $y=y d' d$ and $x=x d d'$. By the cancellative property, $e=d d'$ and $e=d' d$. However, this only holds when $d=d'=e$ by the positive property. Therefore, $x=y$, as desired. (iii) Transitive property. Suppose $x\leq y $ and $y\leq z$. Then there exist $d,d'\in M$ such that $y=x d$ and $z=y d'$. From this, $z=x d d'$, which yields that $x\leq z$ by considering $d d'$. 

    The compatibility of the monoid operation with this order is an immediate consequence of the relation considered.

    Finally, it is a direct consequence of (\ref{eq.leq}) that $e$ is a minimum.
\end{proof}

\subsection{Grothendieck groups} \label{subsection.grothendieck}
Now let us recall some concepts and results concerning the \emph{Grothendieck group}. For a general discussion and for a proof of the following theorem, we refer the reader to \cite{bourbaki}.

\begin{thm}\label{thm_universal_property_grothendieck} Let $(M,\cdot)$ be a commutative monoid.

\begin{itemize}
\item[(1)] Then there exists an abelian group $(K(M),\circledast)$ and a monoid homomorphism $\gamma:M\rightarrow K(M)$ having the following universal property. If $M\rightarrow A$ is a homomorphism into an abelian group $A$, then there exists a unique homomorphism $f_*:K(M)\rightarrow A$ such that $f_*\circ\ \gamma=f$. 

\item[(2)] If the monoid is cancellative, then the monoid homomorphism $\gamma$ is injective and this group can be constructed as $K(M) = (M \times M) / \!\! \sim$,
where the equivalence relation is defined by
\[
(a,b) \sim (c,d) \quad \text{if and only if} \quad  a  d  = b c 
\]
and the group operation $\circledast$ is given by
\[
[(a,b)] \circledast [(c,d)] = [(a c,\, b d)].
\]

\item[(3)]  If the monoid is a unique factorization monoid, then $(K(M),\circledast)$ is free abelian on the set of equivalence classes of 
   irreducible elements of \( M \) modulo units, that is:
   \[
   K(M) \cong \mathbb{Z}^{(I)},
   \]
   where \( I \) indexes the irreducible elements modulo units.

\end{itemize}

\end{thm}

The universal group $K(M)$ is called the \emph{Grothendieck group}. Some immediate consequences of the universal property described in Theorem \ref{thm_universal_property_grothendieck} are:
\begin{cor}\label{cor_extension_homomorfismo_monoides_Grothendieck_groups}
    Let $M,M'$ be monoids. If $f:M\rightarrow M'$ is a monoid homomorphism, then there exists a unique homomorphism $K(f):K(M)\to K(M')$ such that $K(f)\circ \gamma= \gamma'\circ f$. 
\end{cor}

\begin{cor} \label{coro.abelianequal} If $(M,\cdot)$ is not only a monoid but an abelian group, then $K(M)=M$.
    
\end{cor}

\begin{cor} Some particular examples of the Grothendieck group are:

\noindent (a) In the case $(M,\cdot)=(\mathbb N,+)$, the Grothendieck group $K(M)$ can be identified with $\mathbb Z$ under the identification
$$[(a,b)]\sim a-b. $$

\noindent (b) In a similar fashion, for $(M,\cdot)$ equal to $(\mathbb N_{>0},\cdot)$ or $(\mathbb Z\setminus\{0\},\cdot)$, the Grothendieck group $K(M)$ can be identified with $(\mathbb Q_{> 0},\cdot ), (\mathbb Q^\times,\cdot )$, respectively, with respect to the identification
$$[(a,b)]\sim \frac{a}{b}.$$

\end{cor}

We now prove some simple properties that extend from a monoid to its Grothendieck group in a natural manner. However, before that we introduce a useful definition here and in subsequent sections.

\begin{df}\label{def_minimal_representation}
    Let $M$ be a unique factorization cancellative monoid and let $K(M)$ be its Grothendieck group. A representantive $(a,b)\in [(a,b)]$, where $[(a,b)]\in K(M)$, is called a \emph{minimal representation} if the greatest common divisor of $a$ and $b$ is the identity element of $M$.
\end{df}

\begin{lem}\label{lem_minimal_representation_unique}
    Let $M$ be a unique factorization cancellative monoid. Then every element in $K(M)$ admits a minimal representation that is unique.
\end{lem}
\begin{proof}
    Immediate consequence of the unique factorization property and the cancellative property when the greatest common divisor is considered.
\end{proof}

\begin{thm}
    Let $M$ be a monoid with a natural order denoted by $\leq' $. Define a binary relation \( \leq \) on \( K(M) \) by declaring \( [(a,b)] \leq [(c,d)] \) if and only if \( a \leq' c \) and \( b \leq' d \), where \( (a,b) \) and \( (c,d) \) are  minimal representations. Then \( \leq \) is a partial order on \( K(M) \).
\end{thm}
\begin{proof}
    By the uniqueness established in Lemma~\ref{lem_minimal_representation_unique}, the result follows trivially, since we are essentially considering the product order.
\end{proof}

In Section \ref{sec:groth_group_join}, we will recall the construction of $K(M)$ for the specific monoid of simplicial complexes with the join operation, which plays a central role in our exposition.

\subsection{Alexandroff spaces}\label{subsec_alexandroff_spaces}

For a complete introduction on this subject we refer to \cite{barmak2011algebraic,may1966finite}. 

An \emph{Alexandroff space} is a topological space for which the arbitrary intersection of open sets is open. The most important examples are the finite topological spaces.

From now on and without explicit mention, we assume that every Alexandroff space is, additionally, a $T_0$ topological space.

For any Alexandroff space $\mathcal{X}$ and $x\in \mathcal{X}$, we denote by $U_x$ (resp. $F_x$) the intersection of every open (resp. closed) set containing $x$. The relation $x\leq y$ if and only if $U_x\subseteq U_y$ is a partial order on $\mathcal{X}$. Conversely, for every partially ordered set $(\mathcal{X},\leq)$ (or \emph{poset}, for short), the topology on $\mathcal{X}$ generated by considering the basis of the sets $\{U_x:x\in \mathcal{X}\}$, where
$$U_x=\{y\in \mathcal{X}:y\leq x\},$$

\noindent provides a $T_0$ Alexandroff topology for $\mathcal{X}$. These two relations are mutually inverse. Moreover, a map between two Alexandroff spaces is continuous if and only if it is order-preserving. From this, it can be deduced one of the fundamental results in this theory:
\begin{thm}
    The category of Alexandroff spaces and the category of partially ordered sets are isomorphic.
\end{thm}

Therefore, the same object can be seen from two perspectives. In what follows, we will not distinguish between an Alexandroff space or a partially ordered set.

This fact, can be used to describe the notion of homotopy in a simpler manner. An important notion within this context is the notion of \emph{beat point}. A \emph{down (resp. up) beat point} $x$ in an Alexandroff space $\mathcal{X}$ is a point such that $U_x\setminus \{x \}$ has a maximum (resp. $F_x\setminus \{ x\}$ has a minimum).



\begin{prop}\label{prop_contractible} Let $\mathcal{X}$ be an Alexandroff space.

\begin{itemize}

\item[(1)] Removing beat points of $\mathcal{X}$ does not change the homotopy type.

    \item[(2)] If $\mathcal{X}$ contains a maximum or minimum, then it is contractible.

\end{itemize}

\end{prop}

Given an Alexandroff space \( \mathcal{X} \), the \emph{height} of \( \mathcal{X} \), denoted by \( \mathrm{ht}(\mathcal{X}) \), is defined as one less than the maximum number of elements in a chain in \( \mathcal{X} \). For a point \( x \in \mathcal{X} \), the \emph{height} of \( x \), denoted by \( \mathrm{ht}(x) \), is defined as \( \mathrm{ht}(U_x) \).  To conclude, we introduce the following notation: if \( x < y \) in an Alexandroff space \( \mathcal{X} \), and there is no \( z \in \mathcal{X} \) such that \( x < z < y \), then we write \( x \prec y \). If there is no relation between two points $x$ and $y$, we simply write $x\not\sim y$, while $x\sim y$ if $x<y$ or $x>y$.

\subsection{Simplicial complexes}\label{subsec_cw_complexes}

We briefly recall the main objects of study throughout this paper, although we assume the reader is familiar with the basic notions of simplicial complexes. For a comprehensive treatment of these topics, we recommend \cite{munkres1984elements}.

A \emph{simplicial complex} $X$ on a countable vertex set $\{x_1,...,x_n,...\}$ is a collection of subsets of $\{x_1,...,x_n,...\}$ such that for every $\sigma\in X$, if $\tau\subseteq \sigma$, then $\tau \in X$. When the vertex set used to define \( X \) is finite, we say that the simplicial complex \( X \) is finite. Let \( \mathcal{S} \) denote the set of all finite simplicial complexes.

The elements of a simplicial complex $X$ are called \emph{simplices}. The \emph{dimension} of a simplex $\sigma\in X$ is defined as $\dim(\sigma)=\textnormal{Card}(\sigma)-1$ where $\textnormal{Card}(\cdot)$ denotes the cardinality of a set. The \emph{dimension} of $X$, denoted $\dim(X)$, is the the largest dimension among its simplices. Note that in the definition of a simplicial complex, we do not exclude the empty set as an element of $X$; in particular, its dimension is $-1$.

For any positive integer number $n$, consider $[n]$ to be the simplicial complex of dimension 0, containing just the empty set and $n$ singletons corresponding to the vertices. On the other hand, 1-dimensional simplicial complexes are called, in this context, \emph{graphs}.

Given a simplicial complex $X$ of dimension $d$, for every $k=-1,\ldots, d$, the subset $X^k$ consisting in those simplices in $X$ of dimension less than or equal to $k$ is also a simplicial complex called the \emph{$k$-skelleton} of $X$.

Let \( X \) and \( X' \) be simplicial complexes. A map  
$f : X \to X'$ is called a \emph{simplicial map} if it sends vertices of \( X \) to vertices of \( X' \) in such a way that the image of every simplex of \( X \) is a simplex of \( X' \); that is,  
\[
\{v_0, v_1, \dots, v_k\} \in X \implies \{f(v_0), f(v_1), \dots, f(v_k)\} \in X'.
\]

This definition of simplicial complex is purely combinatorial and is commonly referred to in the literature as an \emph{abstract simplicial complex}. There exists a definition of \emph{geometric simplicial complex} which is equivalent in some sense (we omit details). Every abstract simplicial complex $X$ has a geometric counterpart $|X|$ called its \emph{geometric realization}, which is a topological subspace of $\mathbb R^N$, defined up to homeomorphism (we omit details). We say that a simplicial complex is \emph{connected} if so is its geometric realization. Conversely, if a topological space $\mathcal{Y}\subset \mathbb R^N$ is the geometric realization of some simplicial complex $X$, then we say that $X$ is a \emph{triangulation} of $\mathcal{Y}$. In the following, we will not distinguish between simplicial complexes and triangulated topological spaces in $\mathbb R^N$.

Let us recall a useful functor \(\mathcal{X}\) from the category of simplicial complexes to the category of Alexandroff spaces. Given a simplicial complex \(X\), the \emph{face poset} of \(X\), denoted by \(\mathcal{X}(X)\), is the poset whose elements are the simplices of \(X\), ordered by inclusion. Given a simplicial map \(f:X\rightarrow X'\), one has that \( \mathcal{X}(f):\mathcal{X}(X)\rightarrow \mathcal{X}(X')\) defined by \(\mathcal{X}(f)(\sigma)=f(\sigma)\) is a continuous map. 

Let $f:\mathcal{Y}\rightarrow \mathcal{Z}$ be a continuous map between two topological spaces. It is said that $f$ is \emph{a weak homotopy equivalence} if $f$ induces isomorphisms on all homotopy groups. Two spaces $\mathcal{Y}$ and $\mathcal{Z}$ are weak homotopy equivalent if there exists a sequence $f_i:\mathcal{Y}_i\rightarrow \mathcal{Y}_{i+1}$ with $i=0,...,n$ of continuous maps such that $\mathcal{Y}_0=\mathcal{Y}$, $\mathcal{Y}_n=\mathcal{Z}$ and $f_i$ is a weak homotopy equivalence for every $i$. Given a simplicial complex $X$, there exists a natural continuous map $f_X:|X|\rightarrow \mathcal{X}(X)$ that is a weak homotopy equivalence (see \cite{mccord1966singular} for more details).

\subsection{Formal power series}\label{subsec_formal_power_series}

For a complete and comprehensive introduction to this topic, we refer the reader to the book \cite{gan}, the articles \cite{niven,prieto} and the recent preprint \cite{buga}.

 For the rest of this subsection, let \( R \) denote a commutative ring with unit.

 The set of polynomials in one variable \( t \) with coefficients in \( R \) is denoted by \( R[t] \). A \emph{formal power series} with coefficients in \( R \) is an expression of the form
\[
f(t) = \sum_{i=0}^\infty a_i t^i, \quad a_i \in R.
\]
Let \( R[[t]] \) denote the set of all formal power series with coefficients in \( R \). Clearly, \( R[t] \subset R[[t]] \).

It is well known that \( R[[t]] \) is a ring under the usual addition and multiplication of formal power series. An element \( f(t) \in R[[t]] \) is invertible in \( R[[t]] \) if and only if its constant term \( a_0 \) is invertible in \( R \). We denote the group of invertible elements in \( R[[t]] \) by \( R[[t]]^\times \).

A \emph{formal Laurent series} with coefficients in \( R \) is an expression of the form
\[
f(t) = \sum_{n=N}^\infty a_n t^n, \quad a_n \in R, \ N \in \mathbb{Z}.
\]
Let \( R((t)) \) denote the ring of formal Laurent series over \( R \) in the variable \( t \). Similarly, let \( R((1/t)) \) denote the ring of formal Laurent series over \( R \) in \( 1/t \), that is, the set of expressions of the form
\[
f(t) = \sum_{n=-\infty}^N a_n t^n, \quad a_n \in R, \ N \in \mathbb{Z}.
\]

\noindent The \emph{order} of an element \( f(t) \in R((t)) \) (resp. in \( R((1/t)) \)) is the minimal (resp.  maximal) integer \( n \) such that \( a_n \neq 0 \). 

Note that if \( R \) is a field \( \mathbb{K} \), then both \( \mathbb{K}((t)) \) and \( \mathbb{K}((1/t)) \) are fields. Moreover, the sets \( \mathbb{K}((t))^\times = \mathbb{K}((t)) \setminus \{0\} \) and \( \mathbb{K}((1/t))^\times = \mathbb{K}((1/t)) \setminus \{0\} \) form groups under multiplication. 


\section{Grothendieck group of monoids of simplicial complexes with respect to the join operation}\label{sec:groth_group_join}

\subsection{Monoids of simplicial complexes with respect to the join operation are commutative monoids with a natural order}\label{subsec_join_properties}

 The \emph{join} of two simplicial complexes $X$ and $X'$ is the simplicial complex $X*X'$ defined as  $$X*X'= \{\sigma\cup \tau| \sigma\in X, \tau\in X' \}.$$ 

  For a detailed discussion of these definitions and their relationship, see \cite[Chapter 4, Section 2]{matousek2003using}.

The following proposition lists some basic properties satisfied by the join (we refer, again, to \cite{munkres1984elements}).

\begin{prop}\label{prop_propiedades_join_lista}
    Let $X,Y$ and $Z$ be simplicial complexes. Then the following properties hold:
    \begin{enumerate}
        \item Commutativity: $X*Y=Y*X$.
        \item Associativity: $(X*Y)*Z=X*(Y*Z)$.
        \item Identity element: $X*\emptyset=\emptyset*X=X$.
        \item Uniqueness of the unit: $\emptyset$ is the only invertible element.
        \item Cancellation law: if $X*Y=X*Z$ (sim. for $Y*X=Z*X$), then $Y=Z$.
        \item Positiveness: the only invertible element is the identity element $\emptyset$.
        \item $\dim(X*Y)=\dim(X)+\dim(Y)+1$.
    \end{enumerate}
\end{prop}

As an immediate consequence of these properties, we can deduce the following result.
\begin{prop}
    $(\mathcal{S},*)$ is a commutative monoid with a natural order.
\end{prop}

\begin{rem} The set $\mathcal W$ of CW-complexes, endowed with its join operation is a commutative monoid but, unfortunately, it is not cancellative. That hinders to perform an analogue discussion to the one presented here for this more general case. \end{rem}

Let us define $$\mathcal S_{\mathbb N}=\{[n_1]*\ldots*[n_m]|n_i\in \mathbb{N} \ \textnormal{ for every }\ i=1,...,m\ \textnormal{ and some }\  m\in \mathbb{N}\}.$$ Then, we have the following result of immediate proof:

\begin{lem}\label{lem_submonoid} $(\mathcal S_{\mathbb N},*)$ is a submonoid of $(\mathcal{S},*)$ .
\end{lem}

\begin{ex}
    Note that the $n$-dimensional sphere $S^n$ can be triangulated as $[2]*...^{n)}*[2]$, which implies $S^n\in \mathcal{S}_{\mathbb{N}}$. Since $\mathcal{S}_{\mathbb{N}}$ is, by Lemma \ref{lem_submonoid}, a commutative monoid with a natural order, denoted by $\leq$, we have
     $$\emptyset <S^0<S^1<S^2<\cdots S^n<S^{n+1}<\cdots .$$
     This chain contains a homeomorphic copy of every $n$-dimensional sphere. Moreover, it is straightforward to verify that $\textnormal{ht}(S^n)=n+1$.
\end{ex}

Given $X\in \mathcal{S}$, we say that $X$ is a \emph{cone} if there is some $Y\in \mathcal{S}$ such that $X = [1]* Y$ and a \emph{suspension} if there
exists $Z\in \mathcal{S}$ such that $X = [2] * Z$.

By the correspondence recalled in Section \ref{subsec_alexandroff_spaces}, the monoids considered in this section can be interpreted as topological spaces. In the following proposition, we list and state some of their properties.

\begin{prop} The topological spaces $\mathcal{S}$ and $\mathcal{S}_{\mathbb{N}}$ are non-compact, $T_0$, connected, locally connected and contractible. Moreover, $\mathcal{S}_\mathbb{N}$ is an open subset of $\mathcal{S}$.
\end{prop}
\begin{proof}
The only non-trivial part of the proof is the contractibility property. However, this follows from Proposition \ref{prop_contractible}. 
\end{proof}
\begin{rem}
Notice that the minimal closed set containing the point $[1]$ in $\mathcal{S}$ corresponds to all simplicial complexes that are cones, i.e., 
$$F_{[1]}=\{X\in \mathcal{S}| X \textnormal{ is a cone} \}.$$
Similarly, the same situation also holds with the notion of suspension by considering the space $[2]$.
\end{rem}

\subsection{Monoids of simplicial complexes with respect to the join operation are unique factorization monoids}\label{subsec_join_properties2}

Proving that the monoids introduced in Subsection \ref{subsec_join_properties} are, in fact, unique factorization monoids is a bit more subtle.

Let $X\in \mathcal{S}$. We say that $X$ is an \emph{irreducible} simplicial complex if there is no $Y,Z\in \mathcal{S}$, both of them non-empty, such that $X=Z*Y$.

\begin{ex}
    For any positive integer number $n$, the simplicial complex $[n]$ is an irreducible simplicial complex. In fact, we have trivially that any disconnected simplicial complex $X$ is irreducible.
\end{ex}

\begin{thm}\label{thm_caracterizacion_primos_dimension_1}
    Let $X\in \mathcal{S}$ such that $\dim(X)=1$, i.e., $X$ is a graph. Then $X$ is irreducible if and only if $X$ is not a complete bipartite graph.
\end{thm}
\begin{proof} 
    We prove the necessary condition by contraposition. Assume that $X$ is not irreducible of dimension $1$, which implies that there exist non-empty simplicial complexes $Y,Z\in  \mathcal{S}$ such that $X=Y*Z$. By dimensionality, $\dim(Y)=\dim(Z)=0$. Hence, $Y=[k]$ and $Z=[n]$ for some positive integers $k$ and $n$ and, consequently, $X$ is a complete bipartite graph.

    To prove the sufficient condition we argue by contradiction. Suppose $X$ is irreducible and a complete bipartite graph. By the definition of complete bipartite graph $X=[k]*[n]$ for some positive integer numbers, which entails the contradiction.
\end{proof}

The operation of join behaves well with respect to the concept of $k$-skeleton. For every pair of simplicial complexes $X$  and $X'$, the $k$-skeleton of $X*X'$ equals the $k$-skeleton of the join of the $k$-skeletons of $X$ and $X'$. As a consequence:

\begin{prop}\label{prop_1esqueleto_implica_primo} Let $X\in \mathcal{S}$. If the $1$-skeleton of $X$ is not the join of two graphs, then $X$ is irreducible.
\end{prop}
\begin{proof}
    We argue by contradiction. Assume that $X$ is not irreducible. Then there exist non-empty simplicial complexes $Y,Z\in \mathcal{S}$ such that $X=Y*Z$. By considering the $1$-skeletons of $Y$ and $Z$, we obtain a contradiction.
\end{proof}

The converse of Proposition \ref{prop_1esqueleto_implica_primo} does not hold in general as we demonstrate with the following example. 
\begin{ex}\label{ex_irreducible_join_graph_skeleton}
Consider the simplicial complex 
$$X=\{ \{a,b,c \}, \{a,b\},\{a,c\},\{b,c\}, \{d,a\} ,\{d,b\},\{a\},\{b\},\{c\},\{d\},\{e\},\emptyset \},$$

\noindent with vertex set $\{a,b,c,d \}$. We have depicted the simplicial complex $X$ in Figure \ref{fig_irreducible_join_graph_skeleton}. Clearly, $X$ is an irreducible simplicial complex. However, the $1$-skeleton of $X$ is the join of two graphs: $X_1=\{\emptyset,\{a\},\{b\},\{\{a,b \} \}$, and $X_2=\{\emptyset,\{c\},\{d\}\}$.

\begin{figure}[h]
    \centering
    \includegraphics[width=0.3\linewidth]{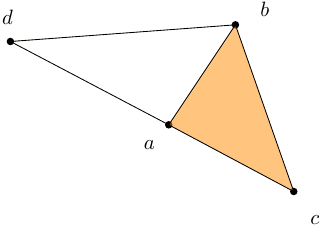}
    \caption{Simplicial complex $X$ from Example \ref{ex_irreducible_join_graph_skeleton}.}
    \label{fig_irreducible_join_graph_skeleton}
\end{figure}
\end{ex}

 Given $X,Y\in \mathcal{S}$, we say that $Y$ \emph{divides} $X$, and we denote it by $Y|X$, if $X=Y*Z$. When $Y$ does not divide $X$, we write $Y\nmid X$.

\begin{rem}\label{rem_cerados}
    Let $X\in \mathcal{S}$. Then the minimal closed set of $X$ corresponds to those simplicial complexes that are divided by $X$, i.e., 
    $$F_X=\{Y\in\mathcal{S}|X|Y \}.$$
\end{rem}

\begin{thm}[Euclid's Lemma]\label{thm_lema_euclides} Let $X\in \mathcal S$. Then $X$ is irreducible if and only if it is \emph{prime}, that is, if $X|(Y*Z)$ for two $Y,Z\in\mathcal S$, then $X|Y$ or $X|Z$. 
\end{thm}

\begin{proof} Notice that prime elements are irreducible in every cancellative monoid, so we omit details.

To prove the converse, by hypothesis, there exists $W\in \mathcal{S}$ such that $X*W=Y*Z$. 

We first prove that $X\subseteq Z$ or $X\subseteq Y$. We argue by contradiction. Assume that $X_1=X\cap Z\neq \emptyset$ and $X_2=X\cap Y\neq \emptyset$. Consider $\sigma_i\in X_i$ for $i=1,2$. Clearly, $\sigma_1\cup \sigma_2$ is a simplex in $Y*Z$ and, consequently, it is a simplex in $X*W$. On the other hand, the simplices of $X*W$ are of the form $\tau\cup \sigma$, where $\tau\in W$ and $\sigma\in X$. Hence, $\sigma_1\cup \sigma_2\in X$ and we deduce that $X=X_1*X_2$, which gives a contradiction with the fact that $X$ is irreducible.

Assume without loss of generality that $X\subseteq Z$. If $X=Z$, then the result holds trivially. Suppose $X\neq Z$ and consider $\tau\in Z\setminus X$. We can deduce that $\tau\in W$ and, generally, that $Z\setminus X\subseteq W$. For every $\tau\in Z\setminus X$ and $\sigma\in X$, we obtain that $\sigma \cup \tau\in X*W=Y*Z$, which implies that $\sigma \cup \tau\in Z$. From here, we deduce that $Z=(Z\setminus X)*X$. By setting $S=Z\setminus X$, we obtain that $X|Z$.

 \end{proof}

\begin{thm}\label{thm_descomposicion_unica_izquierda} The monoid $(\mathcal S,\cdot)$ is a unique factorization monoid, that is, for every $X \in \mathcal{S}$ with $X\neq \emptyset$, there exists a unique decomposition 
$$ X=X_1*\ldots*X_n,$$
where $X_i$ is a prime simplicial complex for $i=1,...,n$. 
\end{thm}
\begin{proof} We argue by induction on the dimension of $X$. 

Suppose the dimension of $X$ is zero. Then $X=[m]$ for some positive integer $m$. Since $[m]$ is a prime simplicial complex we obtain the desired result. 

Suppose the result holds for complexes of dimension $d\geq 0$ and assume that the dimension of $X$ is $d+1$. If $X$ is prime, then the result holds. Suppose $X$ is not prime, which yields that $X$ is not irreducible by Theorem \ref{thm_lema_euclides}. Therefore, there exist non-empty simplicial complexes $Y,Z\in \mathcal{S}$ such that $X=Z*Y$, where both $Z$ and $Y$ have dimension strictly smaller than $d+1$. By the induction hypothesis, 
$$Y=Y_1,\ldots, Y_l,\qquad Z=Z_1*\ldots*Z_n,$$

\noindent where $Y_i,Z_j$ are prime simplicial complexes. Thus,
$$X=Y_1*\ldots*Y_l*Z_1*\ldots*Z_m.$$

\noindent The uniqueness of this decomposition follows trivially from Theorem \ref{thm_lema_euclides}. 


\end{proof}

As an immediate consequence of Theorem \ref{thm_lema_euclides}:
 \begin{cor}
     Let $X,Y\in \mathcal{S}$. If $X*Y$ is a cone, then either $X$ or $Y$ is a cone.
 \end{cor}
\begin{rem}\label{rem_abiertos}
    Given $X\in \mathcal{S}$, the minimal open neighborhood of $X$ corresponds to the set of divisors of $X$. Particularly,
    $$U_X=\{Y\in \mathcal{S}|Y|X\}.$$
    Moreover, if 
    $$X=X_{1_1}*X_{1_2}*\cdots X_{1_{m_1}}*X_{2_1}*\cdots *X_{2_{m_2}}*\cdots *X_{n_1}*\cdots *X_{n_{m_n}}$$ 
    is the decomposition from Theorem \ref{thm_descomposicion_unica_izquierda}, where $X_{i_j}=X_{i_1}$ for every $1\leq j\leq m_j$ and $1\leq i\leq n$, then $$|U_X|=(m_1+1)(m_2+1)\cdots (m_n+1).$$
    This number can be seen as a numerical invariant:  if two simplicial complexes $X$ and $Y$ are isomorphic, then $|U_X|=|U_Y|$.
\end{rem}

Note that the topological space $\mathcal{S}$ can be regarded as a moduli space parametrizing isomorphism classes of finite simplicial complexes.

\subsection{The Grothendieck group of monoids of simplicial complexes with respect to the join operation}

By the results from Subsections \ref{subsec_monoid_grothendieck} and \ref{subsec_join_properties}, we have:

\begin{df} From the monoid $(\mathcal{S},*)$, we define the \emph{Grothendieck group} $(K(\mathcal{S}),+)$, where $K(\mathcal S)$ 
is the set of equivalence classes $[X,Y]$, where, for every $X,Y\in\mathcal S$,
$$[X,Y]=[X',Y']\Leftrightarrow X*Y'=Y*X',\qquad [X,Y]+[X',Y']=[X*X',Y*Y'].$$

\end{df}

To better understand this group, let us highlight the following points:
\begin{itemize}

\item Suppose that $X,Y\in\mathcal S$ have greatest common divisor $Z$ and so $X=X'*Z,Y=Y'*Z$. Then 
$$[X,Y]=\{(X'*R,Y'*R)|R\in\mathcal S\}.$$

\item The identity element of the group is $[\emptyset,\emptyset]$.

\item The inverse of $[X,Y]$ is $[Y,X]$.
\item The map $\gamma:\mathcal S\to K(\mathcal S)$ from Subsection \ref{thm_universal_property_grothendieck} is given by $\gamma(L)=[L,\emptyset]$.

\end{itemize}

As stated before, the group $K(\mathcal S)$ is free abelian and its basis consists of the elements $[X,\emptyset]$ and $[\emptyset,Y]$, where $X$ and $Y$ are prime complexes in $\mathcal{S}$.


We establish the following result that can be deduced from the results of Subsection \ref{subsec_monoid_grothendieck}

\begin{thm}[Fundamental Theorem of Arithmetic]\label{thm_factorizacion_unica_general}
    Let $X,Y\in \mathcal{S}$. Then there exists a unique decomposition 
    $$[X,Y]=\sum_{i=1}^n[X_i,\emptyset]+\sum_{j=1}^m[\emptyset,Y_j], $$
    or equivalently:
    $$X=X_1*\ldots*X_n \textnormal{ and } Y=Y_1*\ldots*Y_m,$$
    where each $X_i$ and $Y_j$ is prime for $i=1,...,n$ and $j=1,...,m$, and such that  $Y_i\nmid X_j$ for all $i,j$. Moreover, if $(Z,W)\in [X,Y]$, then $X\leq Z $ and $Y\leq W$. 
\end{thm}

Let us recall that the representation described in Theorem \ref{thm_factorizacion_unica_general} is called the minimal representation (see Definition \ref{def_minimal_representation}).

\section{Some remarks on the topology of monoids with a natural order and their Grothendieck groups}\label{sec_grothendieck_topology}

\subsection{General remarks}
As discussed in Subsection \ref{subsec_monoid_grothendieck}, for every cancellative and positive monoid $(M,\cdot )$ (this includes $(\mathcal S,*)$), the relation $\leq$, defined in Equation \eqref{eq.leq} is a partial order and, consequently, $(M,\leq)$ is a partially ordered set (or an Alexandroff space). 

Moreover, any such monoid has a minimum (the identity element) and so it is contractible. If such monoid is, additionally, a unique factorization monoid, then we can easily describe the minimal open neighborhood of any element from its prime decomposition. In the case $(M,\cdot)=(\mathcal S,*)$, given an element $X\in \mathcal S$ with prime decomposition $X=X_1*\cdots* X_n$, the minimal open neighborhood $U_{X}$ consists of all possible combinations of joins of the $X_i$ (for $i=1,...,n$). This gives that the space is also locally compact.

Let $ M$ be a unique factorization cancellative monoid and let $K(M)$ be its Grothendieck group. We now introduce a natural topology on $K(M)$, which is fundamentally based on Definition \ref{def_minimal_representation} and Lemma \ref{lem_minimal_representation_unique}. This second asserts that every element $[a,b]\in K(M)$ admits a unique minimal representation.

\begin{thm} Let $(M,\cdot)$ be a unique factorization, positive, cancellative monoid. The group $K(M)$ admits a natural non-discrete Alexandroff topology for which the inversion map is continuous, while the multiplication is not continuous. The corresponding partial order is given by
$$[a,b]\leq [c,d]\qquad\Leftrightarrow a_1\ldots a_n\leq c_1\ldots c_s\text{ and }b_1\ldots b_m\leq d_1\ldots d_t, $$

\noindent where $[a_1\ldots a_n,b_1 \ldots b_m], [c_1\ldots c_s,d_1 \ldots d_t]$ are the minimal representations of $[a,b],[c,d]$.
\end{thm}

\begin{proof}
It is obvious that this relation is a partial order and, therefore, $(K(M),\leq)$ is a partially ordered set.

    The continuity of the inverse map $K(M)\rightarrow K(M)$ defined by $[a,b]\rightarrow [b,a]$ follows trivially: if $[a,b]\leq [c,d]$, then $[b,a]\leq [d,c]$.  However, the operation is not continuous: let $[a,b]\in K(M)$, then $[a,b]\geq [e,e]$, but $[b,a]+[a,b]=[e,e]\ngeq [b,a]=[b,a]+[e,e]$.


\end{proof}

In the following proposition we list some of the topological properties of the above topological space $K(M)$.
\begin{thm}
     Let $(M,\cdot)$ be a non-finite unique factorization, positive, cancellative monoid. Then $K(M)$ is a path connected, non-compact, contractible and locally compact space.
\end{thm}
\begin{proof}
    Since $K(M)$ is non-finite, $K(M)$ is non-compact. The identity element $[e,e]$ is a minimum. This gives that $K(M)$ is contractible. Moreover, for every $[a,b]\in K(M)$, one has that $|U_{[a,b]}|$ is a finite number, which implies that $K(M)$ is locally compact. For any two elements $[a,b],[c,d]\in K(M)$, trivially $[a c,b d]>[a,b],[c,d]$. Hence, $K(M)$ is path connected.
\end{proof}



\subsection{Topology and classes of prime elements}
For the rest of this subsection, $(M,\cdot)$ is a unique factorization, positive, cancellative monoid.

\begin{df}
    We say that $[a,b]\in K(M)$ is the \emph{class of a prime element} if $[a,b]=[a',b']$ where $a'$ or $b'$ is a prime element, while the other is the empty set.
\end{df}

We now provide two characterizations of classes of prime elements in terms of topological or combinatorial aspects.
\begin{prop}
    Let $[a,b]\in K(M)$. Then $[a,b]$ is the class of a prime element if and only if its height is equal to $1$. 
\end{prop}
\begin{proof}
    There are two options when $[a,b]$ is prime: (i) $a\neq e$  is prime and $b=e$, and (ii) $a=e$ and $b\neq e$ is prime. Clearly, we have that the height of $[a,b]$ is  $1$ if (i) or (ii) holds.

    To prove the other implication suppose $[a,b]$ has height $1$. This implies that the maximal chain that can be obtained in $U_{[a,b]}$ has length $1$. Particularly, $[e,e]<[a,b]$ is a chain of length $1$. If $a\neq e$ and $b\neq e$, then we can obtain a chain of length $2$ given by $[e,e]<[a,e]<[a,b]$ which entails a contradiction. Assume without loss of generality that $b=e$. If $a$ is not a prime element, then $a=c d$ for some $c,d\in K(M)$ with $c\neq e\neq d$. However, in that case, we can form the following chain of length $2$: $[e,e]<[c,e]<[a,e]$. This leads to a contradiction and thus proves that, under the given hypothesis, either $a$ or $b$ is a prime element, while the other is the empty set.

\end{proof}

\begin{prop}\label{prop_prime_then_down_beat_point}
    Let $[a,b]\in K(M)$ be the class of a prime element. Then $[a,b]$ is a down beat point.
\end{prop}
\begin{proof} 
    If $[a,b]$ is the class of a prime element, then $a=e$ and $b$ is prime, or $a$ is prime and $b=e$. Without loss of generality assume that $a=e$ and $b$ is a prime. Then $U_{[e,b]}=\{[e,e],[e,b] \}$. It is evident that $U_{[e,b]}\setminus [e,b]$ contains a maximum, which implies that $[e,b]$ is a down beat point.   
\end{proof}
The opposite of Proposition \ref{prop_prime_then_down_beat_point} does not hold as we demonstrate in the following example. 
\begin{ex}
    Let $a\in M$ be prime. The element $[a a,e]$ is a down beat point. However, it is not the class of a prime element. 
\end{ex}

\begin{prop}
    Let $[a,b]\in K(M)$ with minimal representation $[a_1 ...  a_n,b_1 ...  b_m]$. Then $[a,b]$ is a down beat point if and only if $a=e$ and $b_i=b'$ for every $i=1,...,m$ or $b=e$ and $a_i=a'$ for every $i=1,...,n$. 
\end{prop}
\begin{proof}
    Suppose that $[a,b]$ is a down beat point. We study cases. If $n>2$ and $a_i\neq a_j$ for some $1\leq i,j\leq n$, then we may assume that $a_1\neq a_2$. Therefore, $[a_1 a_3  ... a_n,b ]\prec [a,b]$ and $[a_2 ... a_n,b]\prec [a,b]$, which entails a contradiction with the fact that $[a,b]$ is a down beat point. Similarly, we also have the same situation by assuming that $m>2$ and that there are $b_i\neq b_j$ with $1\leq i,j\leq m$. If $n=m=1$ and $a\neq e\neq b$, then $[a,e]\prec [a,b]$ and $[e,b]\prec [a,b]$, which gives the contradiction. Hence, the only possibility is that $a_i=a'$ and $b_j=b'$ for $1\leq i\leq n$ and $1\leq j \leq m$. However, if $a\neq e \neq b$, then $[a,e]\prec [a,b]$ and $[e,b]\prec [a,b]$. This entails a contradiction. Thus, the only option is that $a=e$ and $b_i=b'$ for every $i=1,...,m$ or $b=e$ and $a_i=a'$ for every $i=1,...,n$. 

    Without loss of generality, assume that $a=e$ and $b_i=b'$ for every $i=1,...,m$. Therefore, $[a,b'...^{m-1)} b']$ is a maximum in  $U_{[a,b]}\setminus {[a,b]}$. 
\end{proof}

\subsection{$K(\mathcal S_{\mathbb N})$ as a universal space for the class of finite $T_0$ topological spaces}

To conclude this section, we recall the definition of a \emph{universal space}: a topological space that contains a homeomorphic image of every topological space within a given class. We thus  obtain that $K(\mathcal{S}_{\mathbb N})$ is a universal space for the class of finite $T_0$ topological spaces.

\begin{thm}\label{thm_universal}
    The space $K(\mathcal{S}_{\mathbb{N}})$ (or $K(\mathcal{S})$) is a universal space for the class of finite $T_0$ topological spaces.
\end{thm}
\begin{proof}
    Let $\mathcal{X}$ be a finite $T_0$ topological space such that $\textnormal{ht}(\mathcal{X})=m$. Consider $\mathcal{X}_n=\{x\in \mathcal{X}|\textnormal{ht}(\mathcal{X})=n \}$, where $0\leq n\leq m$. We enumerate the points of $\mathcal{X}$ as follows: 
    \begin{align*}
    \mathcal{X}_0=\{x_0,x_1...,x_{s_0} \}, \mathcal{X}_1=\{ x_{s_0+1}, x_{s_0+2},...,x_{s_{1}}\},...,\mathcal{X}_{m}=\{x_{s_{m-1}+1},x_{s_{m-1}+2},...,x_{s_m} \}.
    \end{align*}
    Then $\mathcal{X}=\bigcup_{i=1}^m \mathcal{X}_i$. Moreover, for every $x\in \mathcal{X}_n$, where $1\leq n\leq m$, set $\mathcal{Y}_n^x=\{y\in \mathcal{X}_{n-1}|y<x\}$. 
    
    The idea of the proof is to define a map $i:\mathcal{X}\rightarrow K(\mathcal{S}_\mathbb{N})$ recursively satisfying that $i$ is a homeomorphism onto its image. Let us start with the base case $n=0$. Define  $i_0:\mathcal{X}_0\rightarrow K(\mathcal{S}_{\mathbb N})$ by $i_0(x_j)=[[j],\emptyset]$ for every $j=1,...,s_0$. Finally, define $i_n:\bigcup_{i=1}^n \mathcal{X}_i\rightarrow K(\mathcal{S_{\mathbb N}})$, for every $1\leq n\leq m$,  by 
   \[i_n(x)=\begin{cases}
        \sum_{y\in \mathcal{Y}_n^x}i_{n-1}(y)+[\emptyset,[l]] \textnormal{ if } x=x_{l}^{n}\in \mathcal{X}_{n},\\
        i_{n-1}(x) \textnormal{ otherwise}.
    \end{cases}
   \]
   Let $i$ denote the map $i_m$ defined recursively above. By construction, $i$ is clearly a bijection onto its image. To conclude, we verify that $x\prec z$ in $\mathcal{X}$ if and only if $i(x)\prec i(z)$ in $i(\mathcal{X})$ because the remaining relations can be derived from transitivity. Suppose $x\prec z$ and $\textnormal{ht}(x)=n$ with $0\leq n<m$, therefore $\textnormal{ht}(z)=n+1$. By construction and hypothesis, $i(x)=i_{n}(x)$, $x\in \mathcal{Y}^z_n$ and 
   $$i(z)=i_{n+1}(z)=\sum_{y\in \mathcal{Y}^z_n} i_{n}(y)=i_{n-1}(x)+ \sum_{y\in \mathcal{Y}^z_n\setminus x} i_{n}(y).$$
   Hence, $i(x)\prec i(z)$. To prove the last implication we argue by contraposition. Suppose $x\not \prec z$. There are three cases: (i) there exists $y$ such that $x<y<z$, (ii) $z\leq x$ and (iii) $x\not \sim z$. If $(i)$ holds, then $i(x)<i(y)<(z)$ and $i(x)\not\prec i(z)$. Similarly, if $(ii)$ holds, then $i(z)\leq i(x)$ and $i(x)\not \prec i(z)$. Finally, assume $(iii)$ holds. Then there exists at least one prime complex $[l]$ in one of the coordinates of $i(\mathcal{X})$ that differs from all prime complexes of $i(z)$.  Similarly, there exists at least one prime complex $[t]$ in $i(z)$ that differs from all prime complexes of  $i(z)$. Therefore, neither $i(x)\leq i(z)$ nor $i(z)\leq i(x)$ holds, and thus $i(x)\not\sim i(z)$.
    
\end{proof}

We provide an example to illustrate the construction of $i:\mathcal X\rightarrow K(\mathcal{S}_\mathbb{N})$ from the proof of Theorem \ref{thm_universal}
\begin{ex}
    Let $\mathcal X=\{x_1,x_2,x_3,x_4 \}$ be the partially ordered set defined by the following parial order $x_1,x_2<x_3,x_4$. Define $i:\mathcal X\rightarrow K(\mathcal{S}_{\mathbb{N}})$ by $i(x_1)=[[1],\emptyset]$, $i(x_2)=[[2,\emptyset]]$, $i(x_3)=[[1],\emptyset]+[[2],\emptyset]+[\emptyset,[3]]=[[1]*[2],[3]]$ and $i(x_4)=[[1],\emptyset]+[[2],\emptyset]+[\emptyset,[4]]=[[1]*[2],[4]]$.
\end{ex}

\begin{rem}
    Despite the universality of $K(\mathcal{S}_{\mathbb{N}})$, this topological space does not provide an alternative description of the \emph{generic poset} $U$. This can be readily verified by observing that some of the properties characterizing $U$ do not hold for $K(\mathcal{S})$. In particular, $K(\mathcal{S})$ fails to satisfy the homogeneity property, which says that every homeomorphism between finite subspaces of $K(\mathcal{S})$ extends to a homeomorphism of the entire space. For further details on the generic poset, we refer the reader to \cite{cameron2010poset}.
\end{rem}

\begin{cor}
    The order complex of $K(\mathcal{S})$ is a universal space for the class of compact ordered complexes.
\end{cor}
\begin{proof}
    It is trivial from Theorem \ref{thm_universal} and the functor $\mathcal{K}:Poset\rightarrow SimComp$ (see \cite{barmak2011algebraic}).
\end{proof}

\section{Notable homomorphisms of groups}\label{sec_homomorphism}

In this section, we construct natural homomorphisms of groups $K(\mathcal{S})\rightarrow K(M)$, where $M$ is a commutative monoid and $K(M)$ denote its Grothendieck group. Specifically, we consider well-known monoid homomorphisms from $(\mathcal{S},*)$ to commutative monoids that encode topological and combinatorial information of simplicial  complexes from $\mathcal{S}$. We then examine the homomorphism $K(\mathcal{S})\rightarrow K(M)$ uniquely determined by Corollary \ref{cor_extension_homomorfismo_monoides_Grothendieck_groups}. This construction enables us to extend classical notions --such as dimension, $f$-vector, and related invariants-- to the elements of $K(\mathcal{S})$. 

\subsection{Dimension}\label{subsec_dimension}

Recall that Proposition \ref{prop_propiedades_join_lista}  states that the join operation can be viewed as a method for constructing, from two simplicial complexes, a new simplicial complex of higher dimension. On the other hand, if $X=Y*Z$ where $X,Y,Z\in \mathcal{S}$, the operation $[X,\emptyset]+[\emptyset,Y]$ results in $[Z,\emptyset]$. This reflects a decrease in dimension. 

Moreover, the dimension of a simplicial complex defines a well-defined map $\dim:\mathcal{S}\rightarrow \mathbb{Z}$, which, however, is not a monoid homomorphism. This is because the map does not send the identity element of $\mathcal{S}$ to the identity element of $\mathbb{Z}$: $\dim(\emptyset)=-1$. More generally, the failure arises from $\dim(X*Y)=\dim(X)+\dim(Y)+1$. Nevertheless, it is straightforward to modify the map $\dim$ to obtain a monoid homomorphism. 
\begin{lem}\label{lem_dim_1_monoid}
    The map $\dim_1:\mathcal{S}\rightarrow \mathbb{N}$ defined by $\dim_1(X)=\dim(X)+1$ is a monoid homomorphism.
\end{lem}
\begin{proof}
    Let $X,Y\in \mathcal{S}$. By applying Proposition \ref{prop_propiedades_join_lista}, we obtain
    $$\dim_1(X*Y)=\dim(X*Y)+1=\dim(X)+\dim(Y)+2=\dim_1(X)+\dim_1(Y),$$
    as desired.
\end{proof}

As stated in Subsection \ref{subsection.grothendieck}, the Grothendieck group $K(\mathbb N)$ can be identified with $(\mathbb Z,+)$. By applying Corollary \ref{cor_extension_homomorfismo_monoides_Grothendieck_groups} to the monoid homomorphis $\dim_1$, we obtain a group homomorphism \begin{equation}\label{eq.homisom} K(\dim_1):K(\mathcal{S})\rightarrow \mathbb{Z}\end{equation}
which is the unique homomorphism satisfying $K(\dim_1)\circ\gamma=\dim_1$. For simplicity, we will henceforth denote $K(\dim_1)$ by $\dim_1$ as well, despite the slight abuse of notation. It is easy to deduce that \[\dim_1([X,Y])=\dim(X)-\dim(Y),\] where $[X,Y]\in K(\mathcal{S})$.

\begin{df}
    We say that the \emph{dimension} of $[X,Y]\in K(\mathcal{S})$ is $\dim_1([X,Y])$.
\end{df}

Note that $\dim_1([X,\emptyset])=\dim(X)+1$, $\dim_1([\emptyset,\emptyset])=0$ and $\dim_1([\emptyset,X])=1-\dim(X)$. Therefore, there are elements from $K(\mathcal{S})$ with negative dimension. 
\begin{ex} 
    Let $X=[2]*[2]$, i.e., the simplicial complex that is a triangulation of the circle $S^1$. We obtain that $\dim_1([X,\emptyset])=2$. Applying a ``desuspension''-- that is,  adding the element $[\emptyset,[2]]$ to $[X,\emptyset]$-- yields the element $[[2],\emptyset]$, which has dimension $1$. Repeating the process, we obtain $[\emptyset,\emptyset]$, of dimension $0$. Finally, applying another  ``desuspension'' to $[\emptyset,\emptyset]$ results in the element $[\emptyset,[2]]$, which has negative dimension $\dim_1([\emptyset,[2]])=-1$. Figure \ref{fig:suspension_desuspension} schematically illustrates this process: applying $[\emptyset,[2]]$ decreases dimension, while applying $[[2],\emptyset]$ increases it.
    \begin{figure}[h!]
        \centering
        \includegraphics[width=0.6\linewidth]{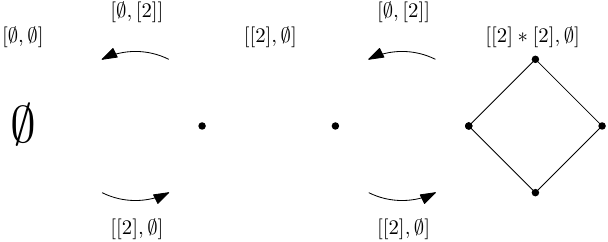}
        \caption{Schematic description of the suspension and ``desuspension''.}
        \label{fig:suspension_desuspension}
    \end{figure}
\end{ex}

\begin{rem}
    The group of integer numbers with the Alexandroff topology provided by the usual order is a topological space. However, the homomorphism $\dim_1:K(\mathcal{S})\rightarrow \mathbb{Z}$ is not a continuous map.
\end{rem}

\subsection{Total reduced Betti number}\label{subsec_total_reduced_betti_number}

Let $\tilde{H}_i(X)$ denote the \emph{reduced $i$-dimensional singular homology group} of a simplicial complex $X$ for $i \geq -1$, with coefficients in a field $\mathbb{K}$. For this and other notions used throughout this section, we refer the reader to \cite{hatcher2000algebraic, munkres1984elements, spanier1981algebraic}. 

In the degenerate case $X = \emptyset$, we assume $\tilde{H}_{-1}(\emptyset) \cong \mathbb{K}$ and $\tilde{H}_i(\emptyset) = 0$ for all $i \geq 0$. Moreover, for $X \neq \emptyset$, we assume $\tilde{H}_{-1}(X) = 0$. For a simplicial complex $X$, let $\tilde{\beta}_i(X) = \operatorname{rank}(\tilde{H}_i(X))$. The number $\tilde{\beta}_i(X)$ is called the \emph{$i$-th reduced Betti number} of $X$. When no confusion arises, we omit the simplicial complex $X$ from the notation. 

Finally, for any simplicial complex $X$, define the \emph{total reduced Betti number} $\tilde{\beta}(X)$ as
\[
\tilde{\beta}(X) = \sum_{i = -1}^{\dim(X)} \tilde{\beta}_i(X). 
\]

The reduced homology of the join of two simplicial complexes $X$ and $Y$ satisfies
\[
\tilde{H}_{i+1}(X * Y) \cong \bigoplus_{j + \ell = i} \left( \tilde{H}_j(X) \otimes \tilde{H}_\ell(Y) \right),
\]
see \cite{bjorner1995topological} for discussions and \cite{milnor1956bundle} for a proof of this result. From this formula, one can deduce a useful property for the total reduced Betti number of the join of two simplicial complexes.
\begin{lem}\label{betti_number_multiplicativo} If $\tilde{\beta}(X)$ and $\tilde{\beta}(Y)$ denote the total reduced Betti numbers of two simplicial complexes $X$ and $Y$, respectively, then
\[
\tilde{\beta}(X * Y) = \tilde{\beta}(X) \cdot \tilde{\beta}(Y).
\]
\end{lem} 
Particularly, the map $\tilde{\beta} : (\mathcal{S},*) \rightarrow (\mathbb{N},\cdot)$ defined by
\[
\tilde{\beta}(X) = \sum_{i = -1}^{\dim(X)} \tilde{\beta}_i(X)
\]
is a monoid homomorphism.  By repeating the argument from Subsection~\ref{subsec_dimension}, we obtain a well-defined group homomorphism
\[
K(\tilde{\beta}) : K(\mathcal{S}) \rightarrow K(\mathbb{N}).
\]

\noindent However, the group $K(\mathbb{N})$, if $\mathbb N$ is endowed with the multiplication, is trivial.

For this reason, we consider the set $\mathcal{S}^\times$ of \emph{non-acyclic simplicial complexes}, that is, consisting in those simplicial complexes with nonzero total reduced Betti number. Clearly, $(\mathcal{S}^\times, *)$ is a submonoid of $(\mathcal{S}, *)$. By repeating the previous argument, we obtain a monoid homomorphim $\tilde{\beta}:K(\mathcal{S}^\times)\rightarrow\mathbb{Q}_{>0}$. Hence, there is a unique homomorphism of groups $K(\tilde{\beta}):K(\mathcal{S}^\times)\rightarrow\mathbb{Q}_{>0}$ such that  $K(\tilde{\beta}) \circ \gamma = \tilde{\beta}$. For simplicity, we will denote $K(\tilde{\beta})$ by $\tilde{\beta}$ as well. It is straightforward to obtain that \[
\tilde{\beta}([X,Y])=\frac{\tilde{\beta(X)}}{\tilde{\beta}(Y)}
\]
where $[X,Y]\in K(\mathcal{S}^\times)$. 

\begin{df}
    We say that that the \emph{total reduced Betti number} of $[X,Y]\in K(\mathcal{S}^\times)$ is $\tilde{\beta}([X,Y])$.
\end{df}

\begin{rem}
     As we did before, let $\mathcal{S}_\mathbb{N}^\times$ denote the set of elements in $\mathcal{S}_\mathbb{N}$ with nonzero total reduced Betti number. This set corresponds to the simplicial complexes obtained as the join of elements of the form $[n]$ where $n\in \mathbb{N}\setminus \{1 \}$.
\end{rem}

\begin{ex}
    Let $X=[2]*\cdots^{n)}*[2]$ and $Y=[2]*\cdots^{m)}*[2]$ for any positive integers $n$ and $m$. Then, $\tilde{\beta}([X,Y])=1$.
    Since $\tilde{\beta}$ is a homotopy invariant, it follows that for any simplicial complexes $X$ and $Y$ that are homotopy equivalent to the $n$-sphere and $m$-sphere, respectively, the total reduced Betti number of $[X,Y]$ is also equal to one.
\end{ex}
\begin{rem}
    The homomorphism $\tilde{\beta}:K(\mathcal{S}^\times)\rightarrow\mathbb{Q}_{>0}$ is clearly surjective. Let us consider the Alexandroff topology provided by the usual order on $\mathbb{Q}_{>0}$. Then the map $\tilde{\beta}:K(\mathcal{S}^\times)\rightarrow \mathbb{Q}_{>0}$ is not continuous.
\end{rem}

\subsection{Extended \(f\)-polynomials}\label{sub_f_polynomials}

Let $X$ be a simplicial complex of dimension $d$ and $(f_{-1},f_0,...,f_d)$ its extended $f$-vector. The \emph{extended $f$-polynomial} of $X$ is defined as
\[ f_X(t)=\sum_{i=-1}^{d} f_it^{i+1}.\]
It is well known that the join operation behaves multiplicatively with respect to extended $f$-polynomials. That is:
\begin{lem}\label{lem_f_polinomio_multiplicativo}
Let $X,Y\in \mathcal{S}$. Then $f_{X*Y}(t)=f_X(t)\cdot f_Y(t)$.
\end{lem}


Note that the extended $f$-polynomial $f_X(t)$ of a simplicial complex $X$ can be viewed as an element of $\mathbb{Q}((1/t))$.

As in other subsections, Lemma \ref{lem_f_polinomio_multiplicativo} and Corollary \ref{cor_extension_homomorfismo_monoides_Grothendieck_groups}  guarantee the existence of a unique homomorphism of groups \[K(f):K(\mathcal{S})\rightarrow \mathbb{Q}((1/t))\] defined by 
\begin{equation}\label{eq.fpolext} K(f)([X,Y])=\frac{f_X(t)}{f_Y(t)}, \end{equation}
such that $K(f)\circ \gamma=f$, where $f:\mathcal{S}\rightarrow \mathbb{Q}((1/t))$ is the monoid homomoprhism defined by $f(X)=f_X(t)$. For simplicity, in the following and commiting an abuse of notation, let $f$ denote $K(f)$.

\begin{rem} For this construction, we have decided that the codomain of $f$ is $\mathbb Q((1/x))$. We have two reasons to do this. The first one is that the formal details are simplified since we can use Corollary \ref{coro.abelianequal}. The second one is that it is coherent with the representation of the extended f-vector of iterative joins as a Riordan matrix (see \cite{chocano} and the example at the end of this section).

\end{rem}

\begin{df}\label{def_extended_f_pol_f_vector_adaptadas}
    We say that the \emph{extended $f$-polynomial} of $[X,Y]\in K(\mathcal{S})$ is $f([X,Y])=\sum_{n=-\infty}^N a_it^i\in \mathbb{Q}((1/x))$, while the \emph{extended $f$-vector} of $[X,Y]$ is $(\ldots ,a_{N+1},a_{N})$.
\end{df}


\begin{prop}\label{prop_order_equals_dimension}
    Let $[X,Y]\in K(\mathcal{S})$. Then the order of $f([X,Y])$ equals $\dim_1([X,Y])$.
\end{prop}
\begin{proof}
    By Subsesction \ref{subsec_dimension}, we know that $\dim_1([X,Y])=\dim(X)-\dim(Y)$. On the other hand, the degree of $f_X(t)$ is $\dim(X)+1$, while the degree of $f_Y(t)$ is $\dim(Y)+1$. Moreover, the order of $\frac{1}{f_Y(t)}$ is $-\dim(Y)-1$. By the product formula in $\mathbb{Q}((1/x))$, we obtain that the order of $\frac{f_X(t)}{f_Y(t)}$ is $\dim(X)-\dim(Y)=\dim_1([X,Y])$, as desired.
\end{proof}

\begin{ex}\label{ex_n_simplex_extended_f_polynomial}
    Let us consider a classical example: the $n$-simplex $\Delta^n=[1]*\cdots^{n+1)}*[1]$. It is well-known that 
    \[ f([\Delta^n,\emptyset])=\prod_{i=1}^{n+1}(1+t)=\sum_{i=0}^{n+1}\binom{n+1}{i}t^i \]
    From this, we can compute the extended $f$-polynomial of the negative simplices $[\emptyset,\Delta^n]$ ($\dim([\emptyset,\Delta^n])=-1-n$):
    \begin{align*}
        f([\emptyset,\Delta^0])&=\frac{1}{1+t}=\sum_{i=-\infty}^{-1} (-1)^{i+1} t^i \\ 
        f([\emptyset,\Delta^1])&=\frac{1}{(1+t)^2}=\sum_{i=-\infty}^{-2} (-1)^{i+1}(i+1)t^i \\
        &\vdots  \\
        f([\emptyset, \Delta^n])&=\frac{1}{(1+t)^{n+1}}=\sum_{i=-\infty}^{-n-1}(-1)^{i+n+1}\binom{-i-1}{-i-n-1}t^i.
    \end{align*}

Moreover, if we store (by rows) the corresponding extended f-vectors of the previous elements in a matrix, then we obtain
$$\mathbf P=\left[\begin{array}{c c c c |c| c c c c} \ddots & & &  & & & \\ 
\hdots & 1& & &  & & & \\
\hdots & -2& 1& &  & & & \\ 
\hdots & 1& -1& 1&  & & & \\
\hline \hdots & 0& 0& 0& 1 & & & \\ 
\hline \hdots & 0& 0& 0& 1 &1 & &\\
\hdots & 0&0 &0 & 1 & 2&1 &\\ 
\hdots & 0&0 &0 & 1 &3 &3 &1\\
& \vdots &\vdots &\vdots  & \vdots &\vdots  &\vdots &\vdots  & \ddots  \end{array}\right].$$

\noindent This is a bi-infinite Riordan matrix \cite{luzon,prieto}. This fact is not merely a curiosity; it has numerous applications. For example, we can exploit the iterative construction pattern of this Riordan matrix $\mathbf{P}$—given by the relation
\[
p_{ij} = p_{i-1,j-1} + p_{i-1,j},
\]
to easily determine any of the coefficients of the extended $f$-polynomials associated with the elements involved.
\end{ex}

The previous example can be easily generalized to simplicial complexes of the form 
\(\Delta_n^{\,k} = [n] * \cdots^{(k)} * [n]\), 
where \(n\) is any positive integer. By storing the \(f\)-vectors of \(\Delta_n^{\,k}\) in a matrix, we obtain a Riordan matrix. These Riordan matrices were studied in \cite{chocano}. Consequently, techniques similar to those mentioned above can be used to determine the coefficients of the extended \(f\)-polynomial associated with the element \([\emptyset, \Delta_n^{\,k}]\).

\begin{rem} By combining the results of Subsection~\ref{subsec_dimension} with those of the present subsection (see Proposition~\ref{prop_order_equals_dimension}), we may naturally speak of simplicial complexes of negative dimension---namely, those of the form $[\emptyset, X]$ with $X \in \mathcal{S}$---and we already know how to compute their associated $f$-vectors, that is, the coefficients of $f([\emptyset, X])$. The idea of negative dimensions has appeared speculatively in the literature in connection with the notion of desuspension; see~\cite{wolcott_mcternan_2012}. In this sense, the present work can be viewed as a partial answer to the question \emph{“Desuspension: what does it lead to?”} However, we are still far from a complete answer, since at present the authors lack, among other things, a geometric interpretation of the negative entries that may occur in the $f$-vectors.
\end{rem}

\subsection{$h$-polynomials}\label{sub_h_polinomios}
As we did in Subsection \ref{sub_f_polynomials}, we can also extend easily the notion of $h$-polynomial and $h$-vector. Let us recall that the \emph{$h$-vector} $(h_0,...,h_d)$ of a complex $X$, whose extended $f$-vector is $(f_{-1},f_0,...,f_{\dim(X)})$, is defined as 
\[h_i =\sum_{k=0}^i (-1)^{i-k}\binom{\dim(X)-k}{\dim(X)-i}f_{k-1},\ \ i=1,...,\dim(X).\]
The \emph{$h$-polynomial} of $X$ is $h_X(t)=\sum_{i=0}^{\dim(X)} h_it^i$. Define $h:\mathcal{S}\rightarrow \mathbb{Q}((1/x))$ by \[h(X)=h_X(t).\] By the properties of the join regarding the $h$-polynomials, $h$ is a well defined monoid homomorphism. By Corollary \ref{cor_extension_homomorfismo_monoides_Grothendieck_groups}, we obtain that $K(h):K(\mathcal{S})\rightarrow\mathbb{Q}((1/x))$ defined as \[ K(h)([X,Y])=\frac{h_X(t)}{h_Y(t)}\]
is the unique homomorphism of groups such that $\gamma\circ K(h)=h$. For simplicity, let $h$ also denote $K(h)$. 

\begin{df}
    We say that the \emph{$h$-polynomial} of $[X,Y]\in K(\mathcal{S})$ is $h([X,Y])$. Moreover, the \emph{$h$-vector} of $[X,Y]$ is given by the coefficients of $h([X,Y])$.
\end{df}     

\subsection{Characteristic polynomial}\label{subsection_characteristic}

Let us introduce the matrices considered in \cite{chocano2025determinant}. Given a simplicial complex $X\in\mathcal S$, consider the face poset $\mathcal{X}(X)=\{x_1,...,x_n \}$ and define the order-complement matrix $\overline{\mathbf Z}_X=(x_{ij})_{i,j=1,...,n}$ by
\[ x_{ij}=\begin{cases}
    0 \quad \textnormal{ if } x_i\leq x_j,\\
    1 \quad \textnormal{otherwise.}
\end{cases}
\]

This matrix encloses all the information of $X$ and some of its invariants are related to invariants of $X$. For example, in \cite{CP}, the authors show that the determinant equals, up to sign, the reduced Euler-characteristic of $X$.

Moreover, for any simplicial complex $X$, let \(p_X(\lambda)\) denote the characteristic polynomial of \(\overline{\mathbf Z}_X\). We say that \(p_X(\lambda)\) is the \emph{charactersitic polynomial} of $X$.  Note that unlike other polynomials considered in this manuscript, \(p_X(\lambda)\) is a topological invariant and unique (Theorem 3.2 in \cite{chocano2025determinant}) in the following sense: \(X\) is homeomorphic to \(Y\) if and only if \( p_X(\lambda)=p_Y(\lambda)\). By the results from \cite{chocano2025determinant}, the map $p:\mathcal{S}\rightarrow \mathbb{Q}((1/\lambda))$ defined by $p(X)=p_{X}(\lambda)$ is a well defined monoid homomorphism. From this and Corollary \ref{cor_extension_homomorfismo_monoides_Grothendieck_groups},  we may conclude the following result:

\begin{lem}
    The map $K(p):K(\mathcal{S})\rightarrow \mathbb{Q}((1/\lambda))$ defined by \[
    K(p)([X,Y])=\frac{p_X(\lambda)}{p_Y(\lambda)}
    \] is a well defined homomorphism of groups and it is the unique satisfying $\gamma\circ K(p_\lambda)=p_\lambda$.
\end{lem}
As we did before, $p_\lambda$ will denote also $K(p_\lambda)$.
\begin{df}
    Let $[X,Y]\in K(\mathcal{S})$. We say that $p([X,Y])$ is the \emph{characteristic polynomial} of $[X,Y]$.
\end{df}

\subsection{Euler characteristic}\label{subsection_euler}

Let $X\in \mathcal{S}$, the \emph{Euler characteristic} of $X$, denoted by $\chi(X)$, is defined as \[\chi(X)=\sum_{i=0}^{\dim(X)} (-1)^{i}f_{i}, \]
where $(f_{-1},f_0,...,f_{\dim(X)})$ is the extended $f$-vector of $X$. The \emph{reduced Euler characteristic} of $X$ is defined as 
\[\tilde{\chi}(X)=\sum_{i=-1}^{\dim(X)} (-1)^{i}f_{i}, \]
or equivalently as $\tilde{\chi}(X)=\chi(X)-1$. 

We have two options to extend the concept of Euler characteristic to $K(\mathcal S)$. Both of them are compatible.

\subsection*{Extension through the concept of genus}

The first option is to use the map $g:\mathcal S\to\mathbb Z$, given by $g(X)=1-\chi(X)=-\tilde{\chi}(X)$. The number $g(X)$ is called the \emph{genus} of $X$ (see \cite[Section 10]{knill2020energy} for this terminology and a discussion about this notion).

It is a well-known result that for any two simplicial complexes $X$ and $Y$, $\chi(X*Y)=\chi(X)+\chi(Y)-\chi(X)\chi(Y)$. This can be derived from the formula of the Euler characteristic of the smash product and the suspension, and also the description of the join in the previous terms (see \cite{hatcher2000algebraic} for details). Hence, it is easy to deduce the following lemma.


\begin{lem}\label{lem_genus}
    Let $X,Y\in \mathcal{S}$. Then $g(X*Y)=g(X)\cdot g(Y)$. 
\end{lem}

Similarly as we did in Subsection \ref{subsec_total_reduced_betti_number}, consider $\mathcal{S}^\times_\chi=\{X\in \mathcal{S}|\chi(X)\neq 1\}$. Note that $\mathcal{S}^\times_\chi$ is again a submonoid, such that $\mathcal S^\times \subset\mathcal{S}^\times_\chi\subset\mathcal S$. By Lemma \ref{lem_genus} and Corollary \ref{cor_extension_homomorfismo_monoides_Grothendieck_groups}, we obtain that \[ K(g):K(\mathcal{S^\times_\chi})\rightarrow \mathbb{Q}^{^\times}\]
defined by \[K(g)[X,Y]=\frac{g(X)}{g(Y)} \]
is the unique homomorphism of groups such that $K(g)\circ \gamma=g$. For simplicity, let $g$ also denote $K(g)$.

\begin{df} \label{def.euler1}
    We say that the \emph{genus} of $[X,Y]\in K(\mathcal{S}^\times_\chi)$ is $g([X,Y])$.
\end{df}

Now, we can define 
\begin{equation}\label{eq.chig}\chi([X,Y])=1-g([X,Y]).\end{equation}

\begin{rem}
    The homomorphism $g:K(\mathcal{S}^\times_\chi)\rightarrow \mathbb{Q}^\times$ is surjective. This can be derived from $g([n])=1-n$ for every $n\geq 2$ and $g(X)=1-A$ where $X$ is homotopy equivalent to the wedge sum of $A+1$ circles with $A\geq 0$.
\end{rem}

\subsection*{Extension through the characteristic polynomial}

For any finite simplicial complex $X$ with extended $f$-vector $(f_{-1},f_0,...,f_d)$, let $n_X=\sum_{i=0}^d f_i$. By the results from \cite{CP}, 
$$p_X(0)=\det(\overline{\mathbf{Z}}_X)=(-1)^{n_X+1}\tilde{\chi}(X).$$

By the results from \cite{chocano}, $\det:\mathcal{S_\chi^\times}\rightarrow \mathbb{Z}$ is a well defined monoid homomorphism. By repeating previous arguments, $\det:K(\mathcal{S}_{\chi}^\times)\rightarrow \mathbb{Q}$ is a homomorphism of groups and it is given by $$\det([X,Y])=\frac{p_X(0)}{p_Y(0)},$$

\begin{df}\label{df_determinante_euler}
    We say that the \emph{determinant} of $[X,Y]\in K(\mathcal{S}_{\chi} ^\times)$ is $\det([X,Y])$.
\end{df}

Moreover, as we did before, we can also define \begin{equation}\label{eq_determinante}
    \chi([X,Y])=1-(-1)^{n_Y-n_X}\det([X,Y]).
\end{equation}

\color{black}
\subsection*{Extension through the extended f-polynomial}


The third option is to use that the Euler characteristic of a simplicial complex $X\in \mathcal{S}$ can be defined via the extended $f$-polynomial $f_X(t)$ just by
$$\chi(X)=1-f_X(-1).$$

This fact motivates the following definition:
\begin{df} \label{def.euler2} The Euler characteristic of $[X,Y]\in K(\mathcal{S})$, denoted by $\chi([X,Y])$, is defined as \begin{equation}\label{eq_f_polinomial}
    \chi([X,Y])= 1-f([X,Y])(-1).
\end{equation}
\end{df}

\begin{rem}\label{rem_lim} Let us remark that, if $ f_Y(-1)=0$, then we interpret that
\[\chi([X,Y])= \lim_{t\to -1}1-f([X,Y])(t).\]

So, we admit $\chi([X,Y])=\pm \infty$.
    
\end{rem}

By developing the expression above using Equation \eqref{eq.fpolext}, we deduce that, for those $Y\in \mathcal S$ such that $\chi(Y)\neq 1$,
\[ \chi([X,Y])=\frac{\chi(X)-\chi(Y)}{1-\chi(Y)}.\]

As an immediate consequence,
$$\forall Y\in \mathcal S,\; \chi([\emptyset, Y])=\frac{\chi(Y)}{\chi(Y)-1}.$$

\subsection*{Coherence between Definitions \ref{def.euler1}, \ref{df_determinante_euler} and \ref{def.euler2}}

Note that every definition for the Euler characteristic $\chi([X,Y])$ (provided by Equations \eqref{eq.chig}, \eqref{eq_determinante}  and \eqref{eq_f_polinomial}) are compatible, for those cases in which $\chi(Y)\neq 1$. 

Moreover, the third option (corresponding to Definition \ref{def.euler2}) also applies for the case $\chi(Y)=1$ when the limit from Remark \ref{rem_lim} exists.

\begin{ex}
    Let us consider the $n$-simplex $\Delta^n$. From Example \ref{ex_n_simplex_extended_f_polynomial}, we conclude that $\chi([\emptyset,\Delta^n])=\infty$ only when $n$ is an even number, since for $n$ odd the limit $\chi([\emptyset,\Delta^n])$ does not exists. Moreover, it is also simple to deduce that $\chi([\emptyset,S^n])$ is equal to $0$ if $n$ is odd and $\chi([\emptyset,S^n])$  is equal to $2$ if $n$ is even.
\end{ex}


\bibliography{bibliografia}
\bibliographystyle{plain}

\end{document}